\documentclass[reqno,10pt,a4paper]{amsart} 

\usepackage{kotex}

\usepackage{amsthm}
\usepackage{amsmath}
\usepackage{amsfonts}
\usepackage{amssymb}

\usepackage{mathrsfs}
\usepackage{stmaryrd} 
\usepackage{mathabx} 
\usepackage[T1]{fontenc}
\usepackage{dsfont}

\usepackage{graphicx}
\usepackage[all,2cell]{xy}
\usepackage{tikz}
\usetikzlibrary{matrix,arrows,decorations.pathmorphing,calc}
\makeatletter
\ifcsname tikz@library@external@loaded \endcsname
\else
\usepackage{tikz-cd}
\fi
\makeatother

\usepackage[colorlinks=true]{hyperref}

\usepackage{url}

\usepackage{setspace}

\usepackage{ifthen}

\usepackage{enumerate}

\usepackage{comment}

\usepackage{bm}

\usepackage{enumitem}

\setlistdepth{9}
\renewlist{itemize}{itemize}{9}
\setlist[itemize,1]{label=\textbullet}
\setlist[itemize,2]{label=--}
\setlist[itemize,3]{label=*}
\setlist[itemize,4]{label=\textbullet}
\setlist[itemize,5]{label=--}
\setlist[itemize,6]{label=*}
\setlist[itemize,7]{label=\textbullet}
\setlist[itemize,8]{label=--}
\setlist[itemize,9]{label=*}

\makeatletter
\@ifclassloaded{book}{
\usepackage[nottoc]{tocbibind}
}{}
\makeatother

\usepackage[normalem]{ulem}
\usepackage{cancel}

\makeatletter
\@ifclassloaded{book}{
  \usepackage{fancyhdr}

  \pagestyle{fancy}

  \fancyhf{}
  \fancyhead[LE,RO]{\bfseries\thepage}
  \fancyhead[LO]{\bfseries\rightmark}
  \fancyhead[RE]{\bfseries\leftmark}

  \addtolength{\headheight}{0.5pt} 
  \fancypagestyle{plain}{%
    \fancyhead{} 
  }
}{}
\makeatother

\makeatletter
\@ifclassloaded{book}{
  \usepackage{minitoc}
}{}
\makeatother

\makeatletter
\@ifclassloaded{book}{
  \usepackage{makeidx}

  \makeindex
}{}

\makeatletter
\@ifclassloaded{book}
{
  \setcounter{tocdepth}{3}
  \setcounter{secnumdepth}{5}
  \setcounter{minitocdepth}{5}
}
{
}

\makeatletter
\@ifclassloaded{book}{
  \newcounter{SaveAppChapter}
  \newcounter{SaveTextChapter}
  \edef\CurrentChapter{\thechapter}

  \newboolean{appchapter}
  \setboolean{appchapter}{false}

  \newcommand{\switchchapter}
  {
    \ifthenelse{
      \boolean{appchapter}
    }
    {
      \setcounter{SaveAppChapter}{\value{chapter}}
      \setboolean{appchapter}{false}
      \setcounter{chapter}{\value{SaveTextChapter}}
      \renewcommand{\thechapter}{\arabic{chapter}}
      \renewcommand{\theHchapter}{\arabic{chapter}}      
    }
    {
      \setcounter{SaveTextChapter}{\value{chapter}}
      \setboolean{appchapter}{true}
      \setcounter{chapter}{\value{SaveAppChapter}}
      \renewcommand{\thechapter}{\Alph{chapter}}
      \renewcommand{\theHchapter}{\Alph{chapter}}      
    }
  }
}{}

\newcounter{SaveAppSection}
\newcounter{SaveTextSection}
\edef\CurrentSection{\thesection}

\newboolean{appsection}
\setboolean{appsection}{false}

\newcommand{\switchsection}
{
  \ifthenelse{
    \boolean{appsection}
  }
  {
    \setcounter{SaveAppSection}{\value{section}}
    \setboolean{appsection}{false}
    \ifthenelse{\equal{\thechapter}{\CurrentChapter}}{}
    {
      \setcounter{SaveTextSection}{0}
      \edef\CurrentChapter{\thechapter}
    }
    \setcounter{section}{\value{SaveTextSection}}
    \renewcommand{\thesection}{\thechapter.\arabic{section}}
    \renewcommand{\theHsection}{\theHchapter.\arabic{section}}
  }
  {
    \setcounter{SaveTextSection}{\value{section}}
    \setboolean{appsection}{true}
    \ifthenelse{\equal{\thechapter}{\CurrentChapter}}{}
    {
      \setcounter{SaveAppSection}{0}
      \edef\CurrentChapter{\thechapter}
    }
    \setcounter{section}{\value{SaveAppSection}}
    \renewcommand{\thesection}{\thechapter.\Alph{section}}
    \renewcommand{\theHsection}{\theHchapter.\Alph{section}}
  }
}

\newcounter{SaveAppSubsection}
\newcounter{SaveTextSubsection}
\edef\CurrentSubsection{\thesubsection}

\newboolean{appsubsection}
\setboolean{appsubsection}{false}

\newcommand{\switchsubsection}
{
  \ifthenelse{
    \boolean{appsubsection}
  }
  {
    \setcounter{SaveAppSubsection}{\value{subsection}}
    \setboolean{appsubsection}{false}
    \ifthenelse{\equal{\thesection}{\CurrentSection}}{}
    {
      \setcounter{SaveTextSubsection}{0}
      \edef\CurrentSection{\thesection}
    }
    \setcounter{subsection}{\value{SaveTextSubsection}}
    \renewcommand{\thesubsection}{\thesection.\arabic{subsection}}
    \renewcommand{\theHsubsection}{\theHsection.\arabic{subsection}}
  }
  {
    \setcounter{SaveTextSubsection}{\value{subsection}}
    \setboolean{appsubsection}{true}
    \ifthenelse{\equal{\thesection}{\CurrentSection}}{}
    {
      \setcounter{SaveAppSubsection}{0}
      \edef\CurrentSection{\thesection}
    }
    \setcounter{subsection}{\value{SaveAppSubsection}}
    \renewcommand{\thesubsection}{\thesection.\Alph{subsection}}
    \renewcommand{\theHsubsection}{\theHsection.\Alph{subsection}}
  }
}

\newcounter{SaveAppSubsubsection}
\newcounter{SaveTextSubsubsection}

\newboolean{appsubsubsection}
\setboolean{appsubsubsection}{false}

\newcommand{\switchsubsubsection}
{
  \ifthenelse{
    \boolean{appsubsubsection}
  }
  {
    \setcounter{SaveAppSubsubsection}{\value{subsubsection}}
    \setboolean{appsubsubsection}{false}
    \ifthenelse{\equal{\thesubsection}{\CurrentSubsection}}{}
    {
      \setcounter{SaveTextSubsubsection}{0}
      \edef\CurrentSubsection{\thesubsection}
    }
    \setcounter{subsubsection}{\value{SaveTextSubsubsection}}
    \renewcommand{\thesubsubsection}{\thesubsection.\arabic{subsubsection}}
    \renewcommand{\theHsubsubsection}{\theHsubsection.\arabic{subsubsection}}
  }
  {
    \setcounter{SaveTextSubsubsection}{\value{subsubsection}}
    \setboolean{appsubsubsection}{true}
    \ifthenelse{\equal{\thesubsection}{\CurrentSubsection}}{}
    {
      \setcounter{SaveAppSubsubsection}{0}
      \edef\CurrentSubsection{\thesubsection}
    }
    \setcounter{subsubsection}{\value{SaveAppSubsubsection}}
    \renewcommand{\thesubsubsection}{\thesubsection.\Alph{subsubsection}}
    \renewcommand{\theHsubsubsection}{\theHsubsection.\Alph{subsubsection}}
  }
}

\makeatletter
\@ifclassloaded{book}{
\numberwithin{section}{chapter}
}{}
\makeatother
\numberwithin{figure}{section}
\numberwithin{table}{section}
\numberwithin{equation}{section}

\newtheorem{theorem}{Theorem}[section]
\newtheorem{lemma}[theorem]{Lemma}
\newtheorem{proposition}[theorem]{Proposition}

\newtheorem{thm-def}[theorem]{Theorem-Definition}

\theoremstyle{definition}
\newtheorem{definition}[theorem]{Definition}

\newtheorem{remark}[theorem]{Remark}

\theoremstyle{remark}

\newcommand{\sSet}{\mathrm{sSet}}

\newdir{ >}{{}*!/-10pt/@{>}}

\makeatletter
\newcommand{\xRightarrow}[2][]{\ext@arrow 0359\Rightarrowfill@{\quad#1\quad}{#2}}
\newcommand{\xLeftarrow}[2][]{\ext@arrow 3228\Leftarrowfill@{\quad#1\quad}{#2}}
\makeatother

\newcommand{\bbP}{\mathbb{P}}
\newcommand{\bbQ}{\mathbb{Q}}
\newcommand{\bbR}{\mathbb{R}}

\newcommand{\bbZ}{\mathbb{Z}}

\renewcommand{\O}{\mathcal{O}}

\newcommand{\uD}{\protect\underline{\mathcal{D}}}
\newcommand{\uE}{\protect\underline{\mathcal{E}}}

\DeclareMathOperator{\BL}{L} 
\DeclareMathOperator{\Bs}{Bs}
\DeclareMathOperator{\codim}{codim}
\DeclareMathOperator{\cosupp}{cosupp}

\DeclareMathOperator{\Div}{Div} 
\DeclareMathOperator{\DivQnn}{Div_{\bbQ\ge 0}} 

\DeclareMathOperator{\mc}{\mathrm{mc}} 
\DeclareMathOperator{\Opz}{O}

\DeclareMathOperator{\sPSh}{sPSh}

\DeclareMathOperator{\supp}{supp}
\newcommand{\Un}[1]{\mathrm{U}(#1)} 

\newcommand{\simQ}{\underset{\bbQ}\sim}

\newcommand{\floor}[1]{\lfloor{#1}\rfloor}

\newcommand{\gangle}[1]{\left\langle{#1}\right\rangle}

\newcommand{\MI}{\mathscr{I}}

\newcommand{\NE}{\mathrm{NE}}

\newcommand{\prfr}[2]{(#1)\(\,\Rightarrow\,\)(#2)}
\newcommand{\prfl}[2]{(#1)\(\,\Leftarrow\,\)(#2)}

\renewcommand{\uD}{\protect\underline{D}}
\renewcommand{\uE}{\protect\underline{E}}
\newcommand{\MM}{\mathrm{M}}

\newcommand{\sMI}{H}
\newcommand{\ZE}{Z}

\begin{document}

\title{Q-divisor and Ampleness}

\author{Seunghun Lee}

\address{Department of Mathematics, Konkuk University,
Kwangjin-Gu Hwayang-dong 1, Seoul 143-701, Korea}

\email{mbrs@konkuk.ac.kr}

\thanks{This research was supported by Basic Science Research Program
  through the National Research Foundation of Korea(NRF) funded by the
  Ministry of Education(no. 2017R1D1A1B03027980)}


\date{\today}


\begin{abstract}
  We give two criteria for a divisor on complex smooth
  projective variety to be ample using the multiplier ideal
  sheaf and the model category.
\end{abstract}

\maketitle

\tableofcontents

\section{Introduction}
\label{sec:introduction-1}

Let \(X\) be a complex smooth projective variety. Let \(L\)
be a divisor on \(X\). The purpose of this note is to give
two criteria for \(L\) to be an ample divisor. The main
object of interest is the set of effective \(\bbQ\)-divisors
associated with \(L\).

\begin{definition}
  \label{def:ample_homotopy:5}
  We let \(\DivQnn(X,L)\) be the set of
  \begin{equation}
    \label{eqn:ample_homotopy:14}
    q E=q\otimes E
  \end{equation}
  in \(\bbQ\otimes\Div(X)\) such that \(q\in \bbQ_{\ge 0}\)
  and \(E\in|mL|\) for some \(m\in \bbZ\).
\end{definition}

The first criterion relies on the multiplier ideal sheaf
(\cite{lazarsfeld-04}).

\begin{definition}
  \label{def:ample_homotopy:7}
  Let \(X\) be a complex smooth projective variety. The
  multiplier ideal sheaf \(\MI(D)\) for an effective
  \(\bbQ\)-divisors \(D\) in \(X\) is defined by
  \begin{equation}
    \label{eqn:ample_homotopy:2}
    \MI(D)=f_{*}\O_{Y}(K_{Y/X}-\floor{D})\subseteq \O_{X}
  \end{equation}
  where \(f:Y\rightarrow X\) is a log resolution for \(D\)
  and \(\floor{D}\) is the round-down of \(D\).
\end{definition}

The multiplier ideal sheaf is an ideal sheaf on \(X\) and
is independent of the choice of \(f\). In a sense it
measures the singularity of \(D\). We refer to
\cite{lazarsfeld-04} for the details.

The following is the key notion.

\begin{definition}
  \label{def:ample_homotopy:1}  
  We say that \(L\) \textbf{satisfies} \textbf{(SoO)} if
  every subvariety \(Z\) of \(X\) is a cosupport of a
  multiplier ideal sheaf \(\MI(D)\) for some
  \(D\in\DivQnn(X,L)\).
  \begin{equation}
    \label{eqn:ample_homotopy:17}
    Z=\cosupp \MI(D)
  \end{equation}
  I.e., for each \(x\in X\), \(x\in Z\)
  iff \(\MI(D)_{x}\subsetneqq \O_{X,x}\).
\end{definition}

The following is our first result. We proved it in
\cite{lee-22} when \(L\) was nef.

\begin{theorem}
  \label{thm:ample_homotopy:1}
  Let \(X\) be a complex smooth projective variety. Let
  \(L\) be a divisor in \(X\). Then the following are
  equivalent.
  \begin{enumerate}
  \item \(L\) is ample.
  \item \(L\) satisfies (SoO).
  \end{enumerate}
\end{theorem}

\begin{remark}
  \label{rem:ample_homotopy:7}
  There is a canonical divisor which is big and nef but not
  ample, and does not satisfies (SoO).  Let \(Y\) be a
  projective threefold in \(\bbP^{4}\) of degree \(d\ge 6\)
  with a \(cA_1\) singularity along a smooth curve \(C\).
  Let
  \(\pi:X\rightarrow Y\) be the crepant resolution of \(Y\)
  at \(C\). We have
  \(K_{X}=\pi^{*}K_{Y}=\pi^{*}\O(d-5)\). \(K_{X}\) is nef
  and big, but not ample. Let \(E\) be the exceptional
  divisor of \(\pi\) in \(X\). Let \(B\) be a curve in \(E\)
  which is mapped to \(C\) as an \(e\)-fold covering with
  \(e\ge 2\).  Then by the Nadel vanishing theorem, there is
  no \(D\in\DivQnn(X,K_{X})\) with \(\cosupp\MI(D)=B\).
\end{remark}

The ampleness is useful partly because it has many
characterizations. In \cite{lee-22}, we showed that one can
detect ample divisors among nef divisors with sheaves in the
context of Grothendieck topoi. Our second result which
relies on Theorem~\ref{thm:ample_homotopy:1} shows that we
can detect the ampleness with simplicial presheaves in the
context of model categories (\cite{quillen-67}).

For the second result, we consider \(\DivQnn(X,L)\) as a
category. A morphism \(D\rightarrow E\) in \(\DivQnn(X,L)\)
is simply an inequality \(D\ge E\). Let \(\Opz(X)\) be the
category of Zariski open sets in \(X\) with an inclusion
\(U\subseteq V\) as a morphism \(U\rightarrow V\).
\(\DivQnn(X,L)\) does not have finite limits. But we can
canonically embed it into a category \(\MM(X,L)\) with
finite limits
(Definition~\ref{def:presheaves_and_Q-divisors:7}).  Given a
functor \(F:\DivQnn(X,L)\rightarrow \Opz(X)\), we denote by
\begin{equation}
  \label{eqn:presheaves_and_Q-divisors:140}
  F^{int}:\MM (X,L)\rightarrow \Opz(X)
\end{equation}
the functor defined by the intersection
(Definition~\ref{def:presheaves_and_Q-divisors:9}).

Quillen introduced the model category as a categorical
framework for homotopy theory (\cite{quillen-67}). We refer
to \cite{hovey-99}, \cite{hirschhorn-03} for the
details. For a small category \(C\), the universal model
category \(\Un C\) is the category of simplicial presheaves
on \(C\) with the Bousfield-Kan model structure
(\cite{dugger-01}).  The universal model categories enjoy a
similar universal property that the free groups do in the
category of groups.  A Quillen adjunction between two model
categories is an adjunction that preserves the model
structures. A Quillen equivalence is a Quillen adjunction
that induces an equivalence between the homotopy categories
of the model categories.

\begin{definition}
  [Definition~\ref{def:homotopical_presentation:1}]
  \label{def:ample_homotopy:6}
  Let \(F:C\rightarrow B\) be a functor between small
  categories \(C\) and \(B\). We say that the functor \(F\)
  \textbf{generates a homotopical presentation of the
  category} \(B\) if the Quillen adjunction
  \begin{equation}
    \label{eqn:ample_homotopy:19}
    F_{*}:\BL_{S_{C}}(\Un C)\rightleftarrows \Un B:F^{*}
  \end{equation}  
  associated with \(F\) is a Quillen equivalence where
  \(\BL_{S_{C}}(\Un C)\) is the left Bousfield localization
  (\cite{hirschhorn-03}) of \(\Un C\) at the set \(S_{C}\)
  of all morphisms in \(C\) mapped to an isomorphism in
  \(B\).
\end{definition}


The following is our second result, which gives a model
theoretic characterization of the ampleness.

\begin{theorem}
  \label{thm:ample_homotopically:1}
  Let \(X\) be a complex smooth projective variety. Let
  \(L\) be a divisor on \(X\). Let
  \(\sMI:\DivQnn(X,L)\rightarrow \Opz(X)\) be the functor
  mapping \(D\) to the support \(\supp\MI(D)\) of the
  multiplier ideal sheaf \(\MI(D)\) for \(D\). Then, the
  following are equivalent.
  \begin{enumerate}
  \item \(L\) is ample.
  \item \(\sMI^{int}\) generates a homotopical presentation of
    \(\Opz(X)\).
  \end{enumerate}
\end{theorem}

When \eqref{eqn:ample_homotopy:19} is a Quillen equivalence,
Dugger (\cite{dugger-01}) calls it a small presentation of
the model category \(\Un B\) by a small category \(C\) and a
set \(S_{C}\) of morphisms.  \(C\) may be considered as a
category of generators and \(S_{C}\) a set of relations.
The category \(\Un B\) is a homotopy colimit completion of
\(B\).


In Section~\ref{sec:ampl-mult-ideal} we prove
Theorem~\ref{thm:ample_homotopy:1}.  In
Section~\ref{sec:proof-theorem-ref} we prepare for the proof
of Theorem~\ref{thm:ample_homotopically:1} and
Section~\ref{sec:example-admiss-subc}. It contains a
necessary and sufficient condition for a functor to generate
a homotopical presentation. It also contains a result that
will be used later in Section~\ref{sec:example-admiss-subc}
to demonstrate how to get a homotopical presentations with a
smaller set of generators. The proof of
Theorem~\ref{thm:ample_homotopically:1} is in
Section~\ref{sec:proof-theorem-ref-2}.

\section{Proof of Theorem~\ref{thm:ample_homotopy:1}}
\label{sec:ampl-mult-ideal}

\prfr 1 2 Let \(Z\) be a subvariety of \(X\). Since \(L\) is
ample, there exists \(n_{0}\) such that
\(\O_{X}(nL)\otimes I_{Z}\) is globally generated for all
\(n\ge n_{0}\). So, if we let \(H=\sum_{i}H_{i}\) for
sufficiently large \(n,k\) and general
\(H_{1},\dots,H_{k}\in |\O_{X}(nL)\otimes I_{Z}|\), we have
\begin{equation}
  \label{eqn:building_Grothendieck_topoi_out_of_multiplier_ideal_sheaves:190}
  Z=\cosupp\MI((1-\varepsilon)H)
\end{equation}
for all sufficiently small \(0<\varepsilon\ll 1\).

\prfl 1 2 We proceed as follows.\footnote{(i) and (ii) were
  proved in \cite{lee-22}.}

\begin{enumerate}
\item[(i)] First we show that if \(L\) satisfies (SoO) then
  \(L\) is big.
\item[(ii)] Next, using (i), we show that if \(L\) satisfies
  (SoO) and is nef then \(L\) is ample.
\item[(iii)] Finally, using (i) and (ii), we show that \(L\)
  satisfies (SoO) then \(L\) is nef.
\end{enumerate}

We denote by \(\simQ\) the \(\bbQ\)-linear equivalence
between \(\bbQ\)-divisors.

\subsection{$L$ is big}
\label{sec:l-big}

Let \(\kappa(L)\) be the iitaka dimension of \(L\).  The number of
subvarieties of \(X\) are uncountable. But if \(\kappa(L)\le 0\) then
there are at most countable numbers of \(\bbQ\)-divisors. So, the
assumption implies that \(\kappa(L)\ge 1\). Suppose that
\(\kappa(L)<\dim X\) holds.  We will derive a contradiction.  Given a
divisor \(D\) in \(X\) we denote by \(\Bs(|D|)\) the set-theoretic
base locus of \(|D|\) and by \(\phi_{|D|}\) the rational map
associated with \(|D|\).

We fix an integer \(m\) that computes \(\kappa(L)\). Then, for each
integer \(k\ge 2\), we have \(\kappa(L)=\dim\phi_{|kmL|}(X)\).  By
Proposition~2.1.21 in~\cite{lazarsfeld-04-0}, we may assume that
\(\Bs(|mL|)=\Bs(|kmL|)\) holds for each integer \(k\ge 1\).  We denote
\(\Bs(|mL|)\) by \(\Bs\) for simplicity.

For each integer \(k\ge 1\), let
\begin{equation}
  \label{eqn:building_Grothendieck_topoi_out_of_an_ample_divisor:116}
  \phi_{k}=\phi_{|kmL|}:X\rightarrow Y_{k}
\end{equation}
be the rational map associated with \(|kmL|\).  Let
\begin{equation}
  \label{eqn:building_Grothendieck_topoi_out_of_an_ample_divisor:117}
  \pi_{k}: X_{k}\rightarrow X
\end{equation}
be a resolution of the indeterminacy of \(\phi_{k}\) with \(X_{k}\)
smooth so that
\begin{equation}
  \label{eqn:building_Grothendieck_topoi_out_of_an_ample_divisor:67}
  (\pi_{k})^{*}|kmL|=|M_{k}|+F_{k}
\end{equation}
where \(|M_{k}|\) is base point free, \(\pi_{k}(F_{k})=\Bs\) and
\begin{equation}
  \label{eqn:building_Grothendieck_topoi_out_of_an_ample_divisor:91}
  \pi_{k}: X_{k}\smallsetminus \supp(F_{k})
  \xrightarrow\cong X\smallsetminus \Bs.
\end{equation}
Let 
\begin{equation}
  \label{eqn:building_Grothendieck_topoi_out_of_an_ample_divisor:118}
  \widetilde\phi_{k}=\widetilde\phi_{|M_{k}|}:
  X_{k}\rightarrow Y_{k}  
\end{equation}
be the morphism associated with \(|M_{k}|\).  Let
\begin{equation}
  \label{eqn:building_Grothendieck_topoi_out_of_an_ample_divisor:119}
  X_{k}\xrightarrow{\psi_{k}} Z_{k} \xrightarrow{\nu_{k}} Y_{k}
\end{equation}
be the stein factorization of \(\widetilde\phi_{k}\).
\begin{equation}
  \label{eqn:building_Grothendieck_topoi_out_of_an_ample_divisor:70}
  \begin{tikzcd}[column sep=large]
    X_{k} \ar[r,"{\psi_{k}}"] \ar[d,"{\pi_{k}}"']
    \ar[dr,"{\widetilde\phi_{k}}" description] & Z_{k} \ar[d,"{\nu_{k}}"]\\
    X \ar[r,"{\phi_{k}}"'] & Y_{k}
  \end{tikzcd}
\end{equation}

For each integer \(k\ge 1\), we choose the following open subvarieties.
\begin{enumerate}
\item \(U_{X,k}\subseteq X\smallsetminus \Bs\). By an abuse of
  notation we use the same notation \(U_{X,k}\) for the corresponding
  isomorphic open subvariety in \(X_{k}\)
  \eqref{eqn:building_Grothendieck_topoi_out_of_an_ample_divisor:91}.
  \(U_{X,k}\subseteq X_{k}\smallsetminus \supp(F_{k})\).
\item \(U_{Y,k}\subseteq Y_{k}\).
\item \(U_{Z,k}\subseteq Z_{k}\).
\end{enumerate}
such that the following two properties hold.
\begin{enumerate}
\item When restricting the diagram
  \eqref{eqn:building_Grothendieck_topoi_out_of_an_ample_divisor:70}
  to these open subvarieties, we have a commutative diagram
  \begin{equation}
    \label{eqn:building_Grothendieck_topoi_out_of_an_ample_divisor:101}
    \begin{tikzcd}[column sep=large]
      U_{X,k} \ar[r,"{\psi_{k}}"] \ar[d,"{\pi_{k}}"',"{\cong}"]
      \ar[dr,"{\widetilde\phi_{k}}" description] & U_{Z,k} \ar[d,"{\nu_{k}}"]\\
      U_{X,k} \ar[r,"{\phi_{k}}"'] & U_{Y,k}
    \end{tikzcd}
  \end{equation}
  of morphisms.  By an abuse of notation, we use the same notations
  for the restrictions of \(\phi_{k}\), \(\widetilde\phi_{k}\),
  \(\psi_{k}\) and \(\nu_{k}\).
\item \(U_{Z,k}\) is a smooth variety and
  \(\psi_{k}:U_{X,k}\rightarrow U_{Z,k}\) is a smooth morphism.
\end{enumerate}

Let \(\Lambda\) be the intersection of two general hypersurfaces in
\(X\).
Then, \(\Lambda\) is irreducible,
\begin{equation}
  \label{eqn:building_Grothendieck_topoi_out_of_an_ample_divisor:94}
  \Lambda\not\subseteq \Bs
\end{equation}
and \(\codim_{U_{X,k}}(\Lambda\cap U_{X,k})=2\)
for all integer \(k\ge 1\).  Since \(\Lambda\) is irreducible, the
assumption implies that there is a \(\bbQ\)-divisor \(D\simQ dL\) for
some \(d\in\bbQ_{\ge 0}\) such that
\begin{equation}
  \label{eqn:building_Grothendieck_topoi_out_of_an_ample_divisor:109}
  \supp(\MI(D))=X\smallsetminus  \Lambda.
\end{equation}
We write \(D=\frac{1}{q}A\)
for a divisor \(A\in|amL|\) and integers \(a,q>0\).  
There is an element \(B\in|(\nu_{a})^{*}(\O_{Y_{a}}(1))|\) satisfying
\begin{equation}
  \label{eqn:building_Grothendieck_topoi_out_of_an_ample_divisor:110}
  (\pi_{a})^{*}(A)-F_{a}=(\psi_{a})^{*}(B).
\end{equation}

Recall that we are assuming that \(1\le \kappa(L)<\dim X\) holds.  We
choose general hypersurfaces in \(Z_{a}\)
\begin{equation}
  \label{eqn:building_Grothendieck_topoi_out_of_an_ample_divisor:105}
  B_{1},\dots,B_{\kappa(L)-1}\quad (\subset Z_{a})
\end{equation}
and general hypersurfaces in \(X_{a}\)
\begin{equation}
  \label{eqn:building_Grothendieck_topoi_out_of_an_ample_divisor:93}
  H_{1},\dots,H_{\dim X-\kappa(L)-1}\quad (\subset X_{a})    
\end{equation}
We define a subvariety \(V\) of \(X_{a}\) by
\begin{equation}
  \label{eqn:building_Grothendieck_topoi_out_of_an_ample_divisor:98}
  V=U_{X,a}\bigcap
  \left\{\cap_{i=1}^{\kappa(L)-1}(\psi_{a})^{*}(B_{i})\right\}
  \bigcap\left\{\cap_{j=1}^{\dim X-\kappa(L)-1} H_{j}\right\}.
\end{equation}
We have
\begin{equation}
  \label{eqn:building_Grothendieck_topoi_out_of_an_ample_divisor:97}
  \MI(D|_{\pi_{a}(V)})=\MI(D)|_{\pi_{a}(V)}
\end{equation}
and \(\dim(\Lambda\cap \pi_{a}(V))=0\).
In particular, for each point \(x\) of
\(\pi_{a}^{-1}(\Lambda)\cap V\),
\begin{equation}
  \label{eqn:building_Grothendieck_topoi_out_of_an_ample_divisor:113}
  \supp(\MI(D|_{\pi_{a}(V)}))=\pi_{a}(V)\smallsetminus\{\pi_{a}(x)\}
\end{equation}
near \(\pi_{a}(x)\).  We define a subvariety \(W\) of \(Z_{a}\) by
\begin{equation}
  \label{eqn:building_Grothendieck_topoi_out_of_an_ample_divisor:103}
  W=U_{Z,a}\cap B_{1}\cap \cdots \cap B_{\kappa(L)-1}.
\end{equation}
By construction, \(V\) is a smooth surface and \(W\) is a smooth
curve.  The restriction of \(\psi_{a}\) to \(V\)
\begin{equation}
  \label{eqn:building_Grothendieck_topoi_out_of_an_ample_divisor:106}
  \psi_{a}|_{V}:V\rightarrow W
\end{equation}
is a smooth morphism by Proposition~III.10.4
in~\cite{hartshorne-77}. In particular, every fiber of
\(\psi_{a}|_{V}\) is a smooth curve.

Let \(x\) be a point of \(\pi_{a}^{-1}(\Lambda)\cap V\).  Because of
\eqref{eqn:building_Grothendieck_topoi_out_of_an_ample_divisor:94},
near \(x\), we have
\begin{equation}
  \label{eqn:building_Grothendieck_topoi_out_of_an_ample_divisor:100}
  (\psi_{a})^{*}(B)=(\pi_{a})^{*}(A)
\end{equation}
and
\begin{equation}
  \label{eqn:building_Grothendieck_topoi_out_of_an_ample_divisor:112}
  (\pi_{a}|_{V})_{*}:\MI(\textstyle\frac{1}{q}(\psi_{a}^{*}(B))|_{V})
  \xrightarrow\cong \MI(D|_{\pi_{a}(V)}).    
\end{equation}
Since the fiber of \(\psi_{a}|_{V}\) over \(\psi_{a}(x)\) is a smooth
curve, near \(x\),
\begin{equation}
  \label{eqn:building_Grothendieck_topoi_out_of_an_ample_divisor:114}
  \MI(\textstyle\frac{1}{q}(\psi_{a}^{*}(B))|_{V})=
  \O_{V}(-\floor{\textstyle\frac{1}{q}(\psi_{a}^{*}(B))|_{V}}).
\end{equation}
But, this is a contradiction to
\eqref{eqn:building_Grothendieck_topoi_out_of_an_ample_divisor:113}.

\subsection{$L$ is ample if $L$ is nef}
\label{sec:l-ample-if}

In this Subsection, we assume that \(L\) is nef.

\subsubsection{$L\cdot C>0$}
\label{sec:l-ample}

We suppose that there is an integral curve \(C\) such that
\(L\cdot C=0\). We will derive a contradiction.

Let \(Z_{L}\) be the set of all integral curves \(C\) in \(X\) such
that \(L\cdot C=0\).  Because \(L\) is big, the Zariski closure of the
union of curves in \(Z_{L}\) is a proper subvariety of \(X\).
\begin{equation}
  \label{eqn:building_Grothendieck_topoi_out_of_an_ample_divisor:124}
  \bigcup_{C\in Z_{L}}C\subsetneqq X
\end{equation}

Now we fix an integral curve \(C\) in \(X\) such that \(L\cdot C=0\).
By Corollary~1.4.41 in~\cite{lazarsfeld-04-0},
there is an integer \(c\) such that
\begin{equation}
  \label{eqn:building_Grothendieck_topoi_out_of_an_ample_divisor:126}
  h^{0}(C,\O_{C}(K_{X}+mL))<c
\end{equation}
holds for each integer \(m\). We choose \(c\) general points
\(\underline x=\{x_{1},\dots,x_{c}\}\) in \(C\).  Let
\(I_{\underline x}\) is the ideal sheaf on \(X\) associated with
\(\underline x\).

Let \(k\) be an integer such that
\(\O_{X}(k)\otimes I_{\underline x}\) is globally generated and the
dimension of the image of the rational map associated with
\(|\O_{X}(k)\otimes I_{\underline x}|\) is \(\dim X\).  Let
\(H_{1},\dots,H_{\dim X-1}\) be general elements of
\(|\O_{X}(k)\otimes I_{\underline x}|\). Let
\begin{equation}
  \label{eqn:building_Grothendieck_topoi_out_of_an_ample_divisor:142}
  B=\bigcap_{k=1}^{\dim X-1}H_{i}.
\end{equation}
Then the curve \(B\) is integral by Theorem~3.3.1
in~\cite{lazarsfeld-04-0}. Because of
\eqref{eqn:building_Grothendieck_topoi_out_of_an_ample_divisor:124},
we may assume that
\begin{equation}
  \label{eqn:building_Grothendieck_topoi_out_of_an_ample_divisor:127}
  L\cdot B>0.
\end{equation}

By assumption, there is an effective divisor \(D\simQ dL\)
for some \(d\in\bbQ_{\ge 0}\) such that
\(\cosupp(\MI(D))= B\). We may assume that the scheme
associated with \(\MI(D)\) is reduced at general
points. Then since \(\dim B=1\) and \(L\) is nef and big we
have surjections
\begin{equation}
  \label{eqn:building_Grothendieck_topoi_out_of_an_ample_divisor:129}
  H^{0}(X,\O_{X}(K_{X}+mL))\rightarrow
  H^{0}(B,\O_{B}(K_{X}+mL))\rightarrow 0
\end{equation}
for all sufficiently large integer \(m\gg 0\) by Theorem 9.4.8 in
\cite{lazarsfeld-04}.
Because of
\eqref{eqn:building_Grothendieck_topoi_out_of_an_ample_divisor:127},
we also have surjections
\begin{equation}
  \label{eqn:building_Grothendieck_topoi_out_of_an_ample_divisor:132}
  H^{0}(B,\O_{B}(K_{X}+mL))\rightarrow H^{0}(\underline
  x,\O_{\underline x}(K_{X}+mL)) \rightarrow 0
\end{equation}
for all sufficiently large integer \(m\gg 0\).

Now,
\eqref{eqn:building_Grothendieck_topoi_out_of_an_ample_divisor:129}
and
\eqref{eqn:building_Grothendieck_topoi_out_of_an_ample_divisor:132}
imply that we have surjections
\begin{equation}
  \label{eqn:building_Grothendieck_topoi_out_of_an_ample_divisor:134}
  H^{0}(C,\O_{C}(K_{X}+mL))\rightarrow H^{0}(\underline
  x,\O_{\underline x}(K_{X}+mL)) \rightarrow 0,
\end{equation}
hence, the inequalities
\begin{equation}
  \label{eqn:building_Grothendieck_topoi_out_of_an_ample_divisor:135}
  h^{0}(C,\O_{C}(K_{X}+mL))\ge c
\end{equation}
for all sufficiently large integer \(m\gg 0\). But it contradicts
\eqref{eqn:building_Grothendieck_topoi_out_of_an_ample_divisor:126}.

\subsubsection{$L\cdot V>0$}
\label{sec:lcdot-v0}

The argument is similar to that of
Section~\ref{sec:l-ample}.  We suppose that \(L\) is not
ample and drive a contradiction.

Let \(d\) be the smallest positive integer such that there is a proper
subvariety \(V\) of dimension \(d\) in \(X\) such that
\begin{equation}
  \label{eqn:building_Grothendieck_topoi_out_of_an_ample_divisor:143}
  L^{d}\cdot V=0.
\end{equation}
By (i) and (ii),
we have \(1<d<\dim X\).  Since \(L^{d}\cdot V=0\) and \(L\) is nef,
\begin{equation}
  \label{eqn:building_Grothendieck_topoi_out_of_an_ample_divisor:144}
  h^{0}(V,\O_{V}(K_{X}+mL))\le a\cdot\frac{m^{d-1}}{(d-1)!}
\end{equation}
for some \(a\in\bbQ_{\ge 0}\) for \(m\gg 0\) by Corollary~1.4.41
in~\cite{lazarsfeld-04-0}.  Let \(b\) be an integer satisfying
\begin{equation}
  \label{eqn:building_Grothendieck_topoi_out_of_an_ample_divisor:145}
  b>a.  
\end{equation}

Let \(H\) be a very ample divisor in \(X\). We choose a general
divisor \(B\in|bH|\) such that \(B\cap V\) is integral. We put
\begin{equation}
  \label{eqn:building_Grothendieck_topoi_out_of_an_ample_divisor:146}
  W=B\cap V.
\end{equation}
Then,
\begin{equation}
  \label{eqn:building_Grothendieck_topoi_out_of_an_ample_divisor:147}
  L^{d-1}W=L^{d-1}\cdot bH\cdot V\ge b
\end{equation}
because \(L^{d-1}\cdot H\cdot V\ge 1\) by the minimality of \(d\).

By assumption, there is an effective divisor \(E\simQ eL\)
for some \(e\in\bbQ_{\ge 0}\) such that
\(\cosupp(\MI(E))=W\).  We denote by \(\ZE\) the scheme
associated with \(\MI(E)\).  By Theorem 9.4.8 in
\cite{lazarsfeld-04}, we have surjections
\begin{equation}
  \label{eqn:building_Grothendieck_topoi_out_of_an_ample_divisor:149}
  H^{0}(X,\O_{X}(K_{X}+mL))\rightarrow
  H^{0}(\ZE,\O_{\ZE}(K_{X}+mL))\rightarrow 0
\end{equation}
for each \(m\) such that \(mL-E\) is nef and big.  \(L|_{\ZE}\) is
ample by our choice of \(d\). So, we have surjections
\begin{equation}
  \label{eqn:building_Grothendieck_topoi_out_of_an_ample_divisor:150}
  H^{0}(\ZE,\O_{\ZE}(K_{X}+mL))\rightarrow
  H^{0}(W,\O_{W}(K_{X}+mL))\rightarrow 0
\end{equation}
for all sufficiently large \(m\gg 0 \).  Now,
\eqref{eqn:building_Grothendieck_topoi_out_of_an_ample_divisor:149} and
\eqref{eqn:building_Grothendieck_topoi_out_of_an_ample_divisor:150}
imply that the canonical map
\begin{equation}
  \label{eqn:building_Grothendieck_topoi_out_of_an_ample_divisor:151}
  H^{0}(V,\O_{V}(K_{X}+mL))
  \rightarrow
  H^{0}(W,\O_{W}(K_{X}+mL))
\end{equation}
is a surjection for all sufficiently large \(m\gg 0\).
Corollary~1.4.41 in~\cite{lazarsfeld-04-0}
and \eqref{eqn:building_Grothendieck_topoi_out_of_an_ample_divisor:147}
imply that
\begin{equation}
  \label{eqn:building_Grothendieck_topoi_out_of_an_ample_divisor:152}
  h^{0}(W,\O_{W}(K_{X}+mL))= c\cdot\frac{m^{d-1}}{(d-1)!}+O(m^{d-2})
\end{equation}
holds for some \(c\ge b\). Then,
\eqref{eqn:building_Grothendieck_topoi_out_of_an_ample_divisor:144}
contradicts
\eqref{eqn:building_Grothendieck_topoi_out_of_an_ample_divisor:145}.

\subsection{$L$ is nef}
\label{sec:l-nef}



We begin with a definition and a series of lemmas. We refer
to Chapter~11 in~\cite{lazarsfeld-04} for the asymptotic
multiplier ideal sheaf. 

\begin{definition}
  \label{def:building_Grothendieck_topoi_out_of_an_ample_divisor:2}
  Let \(X\) be a smooth projective variety. Let \(L\) be an
  divisor of non-negative Iitaka dimension in \(X\).
  \begin{enumerate}
  \item We denote the cosupport of the asymptotic multiplier
    ideal sheaf \(\MI(\| L\|)\) associated to \(L\) by
    \(Z(\| L\|)\).
    \begin{equation}
      \label{eqn:building_Grothendieck_topoi_out_of_an_ample_divisor:189}
      Z(\| L\|)=\cosupp \MI(\| L\|)
    \end{equation}
    I.e. \(Z(\| L\|)\) is the set of all
    points \(p\) of \(X\) such that
    \(\MI(\| L\|)_{p}\subsetneqq\O_{X,p}\).
  \item We denote the union
    \(\bigcup_{m\ge 1} Z(\| mL\|)\) by
    \(Z(L)\).
    \begin{equation}
      \label{eqn:building_Grothendieck_topoi_out_of_an_ample_divisor:216}
      Z(L)=\bigcup_{m\ge 1} Z(\| mL\|)
    \end{equation}
    Note that \(Z(m\|L\|)= Z(\|mL\|)\) and
    \begin{equation}
      \label{eqn:building_Grothendieck_topoi_out_of_an_ample_divisor:228}
      Z(\|mL\|)\subseteq Z(\|(m+1)L\|)
    \end{equation}
    hold for all \(m\ge 1\) by Theorem~11.1.8
    in~\cite{lazarsfeld-04}.
  \end{enumerate}
\end{definition}

\begin{remark}
  \label{rem:building_Grothendieck_topoi_out_of_an_ample_divisor:18}
  Let \(L\) be a big divisor on a smooth projective
  variety. Then \(L\) is nef iff \(\MI(\| mL\|)=\O_{X}\)
  holds for all \(m\ge 1\) by Proposition~11.2.18
  in~\cite{lazarsfeld-04}. The proof of the proposition
  shows that the following
  Lemma~\ref{lem:building_Grothendieck_topoi_out_of_an_ample_divisor:9}
  holds.
\end{remark}

\begin{lemma}
  \label{lem:building_Grothendieck_topoi_out_of_an_ample_divisor:9}
  Let \(X\) be a smooth projective variety. Let \(L\) be a
  big divisor on \(X\). Let \(C\) be an \textbf{integral}
  curve in \(X\). If
  \begin{equation}
    \label{eqn:building_Grothendieck_topoi_out_of_an_ample_divisor:210}
    L\cdot C<0
  \end{equation}
  then 
  \begin{equation}
    \label{eqn:building_Grothendieck_topoi_out_of_an_ample_divisor:187}
    C\subseteq Z(L)
  \end{equation}
  holds.
\end{lemma}
\begin{proof}
  Let \(n=\dim X\). Let \(A\) be a very ample divisor on
  \(X\). Then
  \begin{equation}
    \label{eqn:building_Grothendieck_topoi_out_of_an_ample_divisor:211}
    \O_{X}(K_{X}+nA+mL)\otimes\MI(\| mL\|)
  \end{equation}
  is globally generated for all \(m\ge 1\) by
  Theorem~11.2.13 in~\cite{lazarsfeld-04}.  If
  \eqref{eqn:building_Grothendieck_topoi_out_of_an_ample_divisor:187}
  does not hold then there is an increasing sequence
  \(\{m_{k}\}_{k=1}^{\infty}\) such that
  \(C\subsetneqq Z(\| m_{k}L\|)\).  Thus
  \begin{equation}
    \label{eqn:building_Grothendieck_topoi_out_of_an_ample_divisor:204}
    (K_{X}+nA+m_{k}L)\cdot C\ge 0
  \end{equation}
  for all \(m_{k}\ge 1\).  Then the
  inequalities
  \begin{equation}
    \label{eqn:building_Grothendieck_topoi_out_of_an_ample_divisor:203}
    L\cdot C\ge -\frac{1}{m_{k}}(K_{X}+nA)\cdot C
  \end{equation}
  imply that \(L\cdot C\ge 0\).
\end{proof}


In the following lemma
\eqref{eqn:building_Grothendieck_topoi_out_of_an_ample_divisor:191}
is to be read as "if \(\lambda\) is greater than \(m_{0}\)"
not "\(m_{0}\) is smaller than \(\lambda\)".

\begin{lemma}
  \label{lem:building_Grothendieck_topoi_out_of_an_ample_divisor:11}
  Let \(X\) be a smooth projective variety. Let \(L\) be an
  integral divisor of non-negative Iitaka dimension in
  \(X\). Let \(m_{0}\) be a positive integer.  Let
  \(\lambda\) be a positive rational number. Let \(D\) be an
  effective \(\bbQ\)-divisor in \(X\) satisfying
  \begin{equation}
    \label{eqn:building_Grothendieck_topoi_out_of_an_ample_divisor:190}
    D\simQ \lambda L.
  \end{equation}
  If
  \begin{equation}
    \label{eqn:building_Grothendieck_topoi_out_of_an_ample_divisor:191}
    \lambda\ge m_{0}
  \end{equation}
  then
  \begin{equation}
    \label{eqn:building_Grothendieck_topoi_out_of_an_ample_divisor:205}
    \cosupp\MI(D)\supseteq Z(\| m_{0}L\|)
  \end{equation}
  holds. 
\end{lemma}
\begin{proof}
  There are positive integers \(p,q\) and a divisor
  \(A\in |pL|\) satisfying
  \begin{equation}
    \lambda=\frac{p}{q} \quad \text{and} \quad D=
    \frac{1}{q}A.
  \end{equation}
  By Theorem~11.1.8.(i) in~\cite{lazarsfeld-04}, we have
  \(\MI(m_{0}\| L\|)=\MI(\| m_{0} L\|)\). Let \(k\) be an
  integer that computes \(\MI(m_{0}\| L\|)\). Then \(kp\)
  also computes \(\MI(m_{0}\| L\|)\) by Lemma~11.1.1
  in~\cite{lazarsfeld-04}.  Let \(f:Y\rightarrow X\) be a
  log resolution of \(D\), \(|pL|\) and \(|kp L|\).  By
  \eqref{eqn:building_Grothendieck_topoi_out_of_an_ample_divisor:191},
  we have
  \begin{align}
    \label{eqn:building_Grothendieck_topoi_out_of_an_ample_divisor:193}
    K_{Y/X}-\floor{D}
    &\le K_{Y/X}-\floor{\frac{m_{0}q}{p}D}\\
    &= K_{Y/X}-\floor{\frac{m_{0}}{p}A}\\
    &\le K_{Y/X}-\floor{\frac{m_{0}}{p}F_{p}}\\
    &\le K_{Y/X}-\floor{\frac{m_{0}}{kp}F_{kp}}
  \end{align}
  where \(F_{m}\) is the fixed part of \(f^{*}(|mL|)\).
  Since \(kp\) also computes \(\MI(m_{0}\| L\|)\), we have
  \begin{equation}
    \label{eqn:building_Grothendieck_topoi_out_of_an_ample_divisor:194}
    f_{*}\O_{Y}(K_{Y/X}-\floor{\frac{m_{0}}{kp}F_{kp}})=
    \MI(m_{0}\| L\|).
  \end{equation}
  Thus
  \begin{equation}
    \label{eqn:building_Grothendieck_topoi_out_of_an_ample_divisor:195}
    \cosupp \MI(D)\supseteq \cosupp\MI(m_{0}\| L\|) = Z(\| m_{0}L\|)
  \end{equation}
  holds.
\end{proof}

\begin{remark}
  \label{rem:building_Grothendieck_topoi_out_of_an_ample_divisor:15}
  Lemma~\ref{lem:building_Grothendieck_topoi_out_of_an_ample_divisor:11}
  is meaningful only when \(L\) is not nef or we don't know
  if it is nef because \(L\) is nef iff
  \(Z(\| mL\|)=\emptyset\) holds for all \(m\ge 1\) by
  Proposition~11.2.18 in~\cite{lazarsfeld-04}.
\end{remark}

\begin{lemma}
  \label{lem:building_Grothendieck_topoi_out_of_an_ample_divisor:13}
  Let \(X\) be a smooth projective variety. Let \(L\) be a
  big divisor on \(X\).
  \begin{enumerate}
  \item We assume that 
    \begin{equation}
      \label{eqn:building_Grothendieck_topoi_out_of_an_ample_divisor:206}
      \dim Z(L)\ge 2
    \end{equation}
    or
    \begin{equation}
      \label{eqn:building_Grothendieck_topoi_out_of_an_ample_divisor:208}
      \dim Z(L)=1 \text{ and \(Z(L)\) is reducible}
    \end{equation}
    holds.  Then, there is a positive integer \(m_{0}\) such
    that for every \textbf{integral} curve \(C\) in \(X\)
    and every effective \(\bbQ\)-divisor \(D\) in \(X\)
    satisfying
    \begin{equation}
      \label{eqn:building_Grothendieck_topoi_out_of_an_ample_divisor:224}
      \cosupp\MI(D)=C,
    \end{equation}
    if a rational number \(\lambda\) satisfies
    \begin{equation}
      D\simQ \lambda L
    \end{equation}
    then 
    \begin{equation}
      \label{eqn:building_Grothendieck_topoi_out_of_an_ample_divisor:207}
      \lambda< m_{0}    
    \end{equation}
    holds.
  \item We assume that
    \begin{equation}
      \label{eqn:presheaves_and_Q-divisors:338}
      \dim Z(L)=1 \text{ and \(Z(L)\) is irreducible}
    \end{equation}
    or
    \begin{equation}
      \label{eqn:presheaves_and_Q-divisors:339}
      \dim Z(L)=0
    \end{equation}
    holds. Then there is at most one \textbf{integral} curve
    \(C\) with \(L\cdot C<0\).
  \end{enumerate}
\end{lemma}
\begin{proof}
  (1) Under the assumption of (1), there is a positive
  integer \(m_{0}\) such that
  \begin{equation}
    \dim Z(\|m_{0}L\|)\ge 2
  \end{equation}
  or
  \begin{equation}
    \dim Z(\|m_{0}L\|)=1 \text{ and \(Z(\|m_{0}L\|)\) is reducible}
  \end{equation}
  holds. Since the curve \(C\) is irreducible,
  \(C\supseteq Z(\| m_{0}L\|)\) does not hold.  So,
  Lemma~\ref{lem:building_Grothendieck_topoi_out_of_an_ample_divisor:11}
  implies
  \eqref{eqn:building_Grothendieck_topoi_out_of_an_ample_divisor:207}.

  (2) follows from
  Lemma~\ref{lem:building_Grothendieck_topoi_out_of_an_ample_divisor:9}.
\end{proof}

  Below, we will assume that \(L\) is not nef and derive a
  contradiction.  Now that \(L\) is big by (i)
  Lemma~\ref{lem:building_Grothendieck_topoi_out_of_an_ample_divisor:13}
  is applicable.

  Let \(A\) be a very ample divisor on \(X\) such that
  \(K_{X}+A\) is very ample. We write 
  \begin{equation}
    \label{eqn:presheaves_and_Q-divisors:336}
    H=K_{X}+A.
  \end{equation}
  Let \(r\) be the smallest real number such that \(rH+L\)
  is nef. By our choice of \(r\), we know that \(rH+L\) is
  not ample.  Since \(L\) is not nef, \(r\) is a positive
  real number. If \(r\) is rational then, since \(k(rH+L)\)
  also satisfies (SoO) for a sufficiently large and
  divisible \(k\), \(rH+L\) is ample by (ii),
  a contradiction to the minimality of \(r\).  So, \(r\) is
  irrational.

  For a \(\bbR\)-divisor \(D\) in \(X\), we write
  \begin{equation}
    \label{eqn:building_Grothendieck_topoi_out_of_an_ample_divisor:219}
    D_{<0}=\{\xi\in \overline{\NE}(X)\mid D\cdot \xi<0\}.
  \end{equation}
  For each pair \((p,q)\) of integers with \(p<rq\) there is
  a curve \(C\) with \((pH+qL)\cdot C<0\) because \(rH+L\)
  is nef but not ample.  Then, for all positive
  \(0<\epsilon < r\), the number of numerical classes in
  \(((r-\epsilon)H+L)_{<0}\) represented by an integral
  curve is infinite because \(r\) is irrational. Thus the
  number of numerical classes in \(L_{<0}\) represented by
  an integral curve is also infinite because
  \(((r-\epsilon)H+L)_{<0}\subset L_{<0}\).  So, we rule out
  the cases \eqref{eqn:presheaves_and_Q-divisors:338} and
  \eqref{eqn:presheaves_and_Q-divisors:339} in
  Lemma~\ref{lem:building_Grothendieck_topoi_out_of_an_ample_divisor:13}(2).
  Below, we will show that the number of numerical classes
  in \(\overline\NE(X)\) represented by an integral curve is
  finite under
  \eqref{eqn:building_Grothendieck_topoi_out_of_an_ample_divisor:206}
  or
  \eqref{eqn:building_Grothendieck_topoi_out_of_an_ample_divisor:208}
  in
  Lemma~\ref{lem:building_Grothendieck_topoi_out_of_an_ample_divisor:13}(1)
  when \(L\) satisfies (SoO). This will provide a
  contradiction.
  
  Because \(L\) satisfies (SoO), for each integral curve
  \(C\) in \(X\) there are an effective \(\bbQ\)-divisor
  \(D_{C}\) and a positive rational number
  \(\lambda_{C}\in\bbQ_{>0}\) satisfying
  \begin{equation}
    D_{C}\simQ\lambda_{C} L\quad\text{and}\quad
    \cosupp\MI(D_{C})=C.
  \end{equation}
  We fix one \(D_{C}\) for each \(C\).  Let \(m_{0}\) be the
  integer in (1) of
  Lemma~\ref{lem:building_Grothendieck_topoi_out_of_an_ample_divisor:13}.
  Since \(H\) is ample there are positive integers \(p,q\)
  such that
  \begin{enumerate}
  \item[(A)] \(q\ge m_{0}\) 
  \item[(B)] \(pH+qL\) is ample.
  \end{enumerate}  
  Since \(\lambda_{C}<m_{0}\) by
  \eqref{eqn:building_Grothendieck_topoi_out_of_an_ample_divisor:207},
  (A) and (B) imply that
  \begin{enumerate}
  \item[(C)] \(pH+qL-D_{C}\) is ample for all \(C\).
  \end{enumerate}

  Now, we fix an integral curve \(C\) in \(X\). Let
  \(A_{0}\) be a general member of \(|A|\).  Then by (C) and
  the Nadel vanishing theorem
  \begin{align}
    0&=h^{1}(A_{0},\O_{A_{0}}(K_{A_{0}}+pH+qL)\otimes\MI(D_{C}|_{A_{0}}))\\
     &=h^{1}(A_{0},\O_{A_{0}}((p+1)H+qL)\otimes\MI(D_{C}|_{A_{0}}))
  \end{align}
  holds.
  Thus we have a surjection
  \begin{align}
    &\,H^{0}(A_{0},\O_{A_{0}}((p+1)H+qL))\\
    \rightarrow
    &\,H^{0}(A_{0}\cap C, \O_{A_{0}\cap C}((p+1)H+qL))
  \end{align}
  and
  \begin{align}
    &\, h^{0}(A_{0},\O_{A_{0}}((p+1)H+qL))\\
    \ge
    &\, h^{0}(A_{0}\cap C,\O_{A_{0}\cap C}((p+1)H+qL))\\
    \ge
    &\,\deg_{A}C
  \end{align}
  holds.  By (B), we have a surjection
  \begin{equation}
    \label{eqn:building_Grothendieck_topoi_out_of_an_ample_divisor:218}
    H^{0}(X,\O_{X}((p+1)H+qL))
    \rightarrow
    H^{0}(A_{0}, \O_{A_{0}}((p+1)H+qL)).
  \end{equation}
  Thus
  \begin{equation}
    \label{eqn:building_Grothendieck_topoi_out_of_an_ample_divisor:209}
    h^{0}(X,\O_{X}((p+1)H+qL))\ge \deg_{A}C.
  \end{equation}
  Therefore, the set of numerical classes in \(\overline\NE(X)\)
  represented by an integral curve is finite.

\section{Homotopical presentations of small categories}
\label{sec:proof-theorem-ref}

This section states two results in \cite{lee-24-a}
concerning the homotopical presentation of small categories,
which prepare for the proof of
Theorem~\ref{thm:ample_homotopically:1} and
Section~\ref{sec:example-admiss-subc}.

First, we recall a notion from \cite{lee-24-a}.

\begin{definition}
  \label{def:homotopical_presentation:1}
  Let \(F:C\rightarrow B\) be a functor between small
  categories \(C\) and \(B\). Let \(S\) be a set of
  morphisms in \(C\) mapped to an isomorphism in \(B\).
  \begin{itemize}
  \item We say that the functor \(F\) \textbf{generates a
      homotopical presentation of the category \(B\) with a
      generator \(C\) and a relation \(S\)} if the Quillen
    adjunction
    \begin{equation}
      \label{eqn:presheaves_and_Q-divisors:135}
      F_{*}:\BL_{S}(\Un C)\rightleftarrows \Un B:F^{*}
    \end{equation}  
    associated with \(F\) is a Quillen equivalence where
    \(\BL_{S}(\Un C)\) is the left Bousfield localization
    (\cite{hirschhorn-03}) of \(\Un C\) at \(S\).
  \item When the set \(S\) in
    \eqref{eqn:presheaves_and_Q-divisors:135} is the set of
    all morphisms in \(C\) mapped to an isomorphism in
    \(B\), we simply say that the functor \(F\)
    \textbf{generates a homotopical presentation of the
      category \(B\)}.
  \end{itemize}
\end{definition}

Given a small category \(A\), we denote by \(\sPSh(A)\) the
category of simplicial presheaves on \(A\).  The universal
model category \(\Un A\) on \(A\) (\cite{dugger-01}) is the
model category of \(\sPSh(A)\) and the Bousfield-Kan model
structure; a fibration is an object-wise fibration in
\(\sSet\) and a weak equivalence is an object-wise weak
equivalence in \(\sSet\).  Given a functor
\(F:C\rightarrow B\) between small categories, there is an
adjunction
\begin{equation}
  \label{eqn:presheaves_and_Q-divisors:168}
  F_{*}:\sPSh(C)\rightleftarrows \sPSh(B):F^{*}
\end{equation}
obtained from \(F\).\footnote{We follow the notations in
  \cite{artin-62}. Instead of
  \eqref{eqn:presheaves_and_Q-divisors:168},
  \(F^{\dag}:\sPSh(C)\rightleftarrows \sPSh(B):F_{*}\) is used in
  Definition~2.3.1 in \cite{kashiwara-shapira-06}. We do not need the
  right Kan extension functor \(F^{\ddag}\).}  It is the underlying
adjunction of the Quillen adjunction
\begin{equation}
  \label{eqn:presheaves_and_Q-divisors:136}
  F_{*}:\Un C\rightleftarrows \Un B:F^{*}.
\end{equation}
The Quillen adjunction
\eqref{eqn:presheaves_and_Q-divisors:135} is obtained from
\eqref{eqn:presheaves_and_Q-divisors:136} by localizing
\(\Un C\) at \(S\).  Even if we localize with respect to all
morphisms in \(C\) mapped to an isomorphism in \(B\), we
need to impose a connectivity of fibers of \(F\) for
\eqref{eqn:presheaves_and_Q-divisors:135} to be a Quillen
equivalence. For this, we need the following notions.

\begin{definition}
  \label{def:homotopical_presentation:2}
  Let \(C,B\) be categories. Let \(F:C\rightarrow B\) be a
  functor.  Let \(d\) be an object of \(C\).
  \begin{enumerate}
  \item We denote by \(C_{d}\) the full subcategory of the
    comma category \((C\downarrow d)\) such that for each
    object \(\gangle{e,\phi:e\rightarrow d}\) of
    \((C\downarrow d)\), \(\gangle{e,\phi}\) is an object of
    \(C_{d}\) iff \(F\phi\) is an isomorphism in \(B\).
  \item For a subcategory \(A_{d}\) of \(C_{d}\), we say
    that \(F\) \textbf{has \(A_{d}\)-lifting property} if
    \begin{itemize}
    \item for each \(e\in C\) and a morphism
      \(\delta:Fd\rightarrow Fe\), there is an object
      \(\gangle{e',\phi:e'\rightarrow d}\) of \(A_{d}\) and
      a morphisms \(\psi:e'\rightarrow e\) in \(C\) such
      that \(\delta=F\psi \cdot (F\phi)^{-1}\) holds.
    \end{itemize}
  \item For a family of subcategories \(A_{d}\) of \(C_{d}\)
    for \(d\in C\), we denote by \(S_{A}\) the set of
    morphisms \(\phi:e\rightarrow d\) in \(C\) such that
    \(\gangle{e,\phi}\in A_{d}\).  In
    particular, \(S_{C}\) is the set of all morphisms in
    \(C\) mapped to an isomorphism in \(B\).
  \end{enumerate}
\end{definition}

The following result is well-known from the infinity
categorical perspective (\cite{cisinski-19},
\cite{lurie-09}). For a direct proof without using infinity
categories, we refer to \cite{lee-24-a}.

\begin{theorem}[Corollary~1.2 in~\cite{lee-24-a}]
  \label{thm:presheaves_and_Q-divisors:6}
  Let \(B,C\) be small categories. Let \(F:C\rightarrow B\) be a
  functor. We assume that \(C\) has finite limits and \(F\) preserves
  finite limits.  Then the following are equivalent.
  \begin{enumerate}
  \item\label{item:20} The following two properties hold.
    \begin{enumerate}
    \item \(F\) is essentially surjective on objects.
    \item \(F\) has \(C_{d}\)-lifting property for all
      \(d\in C\).
    \end{enumerate}
  \item\label{item:21} \(F\) generates a homotopical presentation of
    \(B\).
  \end{enumerate}
\end{theorem}

The implication \prfr 1 2 in
Theorem~\ref{thm:presheaves_and_Q-divisors:6} can be
generalized with the following notion.

\begin{definition}
  \label{def:ample_homotopy:3}
  Let \(C,B\) be categories. Let \(F:C\rightarrow B\) be a
  functor.  Let \(d\) be an object of \(C\).  We say that a
  subcategory \(A_{d}\) of \(C_{d}\) is
  \textbf{\(F\)-admissible} if \(A_{d}\) is cofiltered and
  \(F\) has \(A_{d}\)-lifting property.
\end{definition}

For the proof of the following result, we also refer to
\cite{lee-24-a}.

\begin{theorem}[Theorem~1.1 in~\cite{lee-24-a}]
  \label{thm:presheaves_and_Q-divisors:1}
  Let \(B,C\) be small categories. Let \(F:C\rightarrow B\) be a
  functor. We assume that \(C\) has finite limits and \(F\) preserves
  finite limits.  We assume that the following hold.
  \begin{enumerate}
  \item \(F\) is essentially surjective on objects. 
  \item For each \(d\in C\), there is a \(F\)-admissible
    subcategory \(A_{d}\) of \(C_{d}\).
  \end{enumerate}
  Then \(F\) generates a homotopical presentation of the
  category \(B\) with a generator \(C\) and a relation
  \(S_{A}\).
\end{theorem}


\section{Proof of Theorem~\ref{thm:ample_homotopically:1}}
\label{sec:proof-theorem-ref-2}

We first recall from \cite{lee-22} the finite limit completion
\begin{equation}
  \label{eqn:presheaves_and_Q-divisors:49}
  \MM (X,L)  
\end{equation}
of \(\DivQnn(X,L)\).\footnote{Beside the notational
  difference from \(M_{X,L}\) in~\cite{lee-22}, there is one
  more minor change. We do not assume that objects of
  \(\MM(X,L)\) to be reduced. Then \(M_{X,L}\) is a skeleton
  of \(\MM(X,L)\).}

\begin{definition}
  \label{def:presheaves_and_Q-divisors:8}
  Let \(X\) be a complex smooth projective variety. Let
  \(L\) be a divisor on \(X\).  Let
  \(\uD=\{D_{1},\dots,D_{n}\}\) and
  \(\uE=\{E_{1},\dots,E_{m}\}\) be two finite sets of
  objects of \(\DivQnn(X,L)\).  We write
  \begin{equation}
    \label{eqn:ample_homotopy:1}
    \uD\ge\uE
  \end{equation}
  if for each \(E_{i}\) there is \(D_{j}\) such that
  \(D_{j}\ge E_{i}\).
\end{definition}

\begin{definition}
  \label{def:presheaves_and_Q-divisors:7}
  Let \(X\) be a complex smooth projective variety. Let \(L\) be a
  divisor on \(X\). We define a category \(\MM (X,L)\) as follows.
  \begin{enumerate}
  \item An object of \(\MM (X,L)\) is a finite set
    \begin{equation}
      \label{eqn:ample_homotopy:16}
      \uD=\{D_{1},\dots,D_{n}\}
    \end{equation}
     of objects of
    \(\DivQnn(X,L)\).
  \item A morphism \(\uD\rightarrow \uE\) in \(\MM (X,L)\) is
    an inequalities \(\uD\ge \uE\) in
    \eqref{eqn:ample_homotopy:1}.
  \end{enumerate}
\end{definition}

\begin{lemma}
  \label{lem:presheaves_and_Q-divisors:14}
  Let \(X\) be a complex smooth projective variety. Let \(L\) be a
  divisor on \(X\). Then the category \(\MM (X,L)\) has finite
  limits. In particular
  \begin{equation}
    \label{eqn:presheaves_and_Q-divisors:13}
    \uD\times\uE=\uD\cup\uE
  \end{equation}
  holds for all objects \(\uD,\uE\) of \(\MM (X,L)\).
\end{lemma}

\begin{remark}
  \label{rem:ample_homotopy:3}
  We consider \(\DivQnn(X,L)\) as a subcategory of
  \(\MM(X,L)\) by mapping \(D\) to \(\{D\}\) for each
  \(D\in \DivQnn(X,L)\).
\end{remark}

\begin{definition}
  \label{def:ample_homotopy:4}
  Let \(X\) be a complex smooth projective variety.  We let
  \(\Opz(X)\) be the category of Zariski open sets in \(X\)
  with a inclusion \(U\subseteq V\) as a morphism
  \(U\rightarrow V\).
\end{definition}

\begin{remark}
  \label{rem:ample_homotopy:4}
  In Definition~\ref{def:ample_homotopy:4}, the limits in
  \(\Opz(X)\) will always be the intersections.
\end{remark}

Next, we recall the extension
\begin{equation}
  \label{eqn:presheaves_and_Q-divisors:37}
  F^{int}:\MM(X,L)\rightarrow \Opz(X)
\end{equation}
of a functor \(F:\DivQnn(X,L)\rightarrow \Opz(X)\) by
intersection.

\begin{definition}
  \label{def:presheaves_and_Q-divisors:9}
  Let \(X\) be a complex smooth projective variety. Let \(L\) be a
  divisor on \(X\). Let \(F:\DivQnn(X,L)\rightarrow \Opz(X)\) be a
  functor. We define a functor
  \begin{equation}
    \label{eqn:presheaves_and_Q-divisors:80}
    F^{int}:\MM (X,L)\rightarrow \Opz(X)
  \end{equation}
  as follows.
  \begin{enumerate}
  \item For each object \(\uD=\{D_{1},\dots,D_{n}\}\) of \(\MM (X,L)\),
    \begin{equation}
      \label{eqn:presheaves_and_Q-divisors:82}
      F^{int}\uD=\cap_{i=1}^{n}FD_{i}.
    \end{equation}
  \item For each morphism \(\phi:\uD\rightarrow \uE\) in
    \(\MM (X,L)\), \(F^{int}\phi\) is the canonical morphism
    in \(\Opz(X)\) induced by \(F\).
  \end{enumerate}
  We call \(F^{int}\) \textbf{the extension of} \(F\)
  \textbf{by intersection}.
\end{definition}

\begin{lemma}
  \label{lem:presheaves_and_Q-divisors:15}
  Let \(X\) be a complex smooth projective variety. Let
  \(L\) be a divisor on \(X\). Let
  \(F:\DivQnn(X,L)\rightarrow \Opz(X)\) be a functor. Then
  the extension \(F^{int}\) of \(F\) by intersection
  preserves finite limits.
\end{lemma}
\begin{proof}
  The limits in \(\Opz(X)\) are intersections. So, it
  follows from \eqref{eqn:presheaves_and_Q-divisors:13} and
  \eqref{eqn:presheaves_and_Q-divisors:82}.
\end{proof}

\begin{lemma}
  \label{lem:presheaves_and_Q-divisors:16}
  Let \(X\) be a complex smooth projective variety. Let
  \(L\) be a divisor on \(X\). Let
  \(F:\DivQnn(X,L)\rightarrow \Opz(X)\) be a functor. Then
  the extension \(F^{int}\) of \(F\) by intersection has
  \(\MM(X,L)_{\uD}\)-lifting property for all
  \(\uD\in\MM(X,L)\).
\end{lemma}
\begin{proof}
  Let \(\uD,\uE\in \MM (X,L)\).  Let
  \(\delta:F^{int}(\uD)\rightarrow F^{int}(\uE)\) be a
  morphism in \(\Opz(X)\). Let
  \(\phi:\uD\times \uE\rightarrow \uD\) and
  \(\psi:\uD\times \uE\rightarrow \uE\) be the
  projections. Then \(F^{int}(\phi)\) is the identity by
  \eqref{eqn:presheaves_and_Q-divisors:13} and
  \eqref{eqn:presheaves_and_Q-divisors:82}.  So, \(F^{int}\)
  satisfies (Lift).
\end{proof}

Finally, we give the proof of
Theorem~\ref{thm:ample_homotopically:1}.

\begin{proof}[Proof of
  Theorem~\ref{thm:ample_homotopically:1}]

  The functor \(\sMI:\DivQnn(X,L)\rightarrow \Opz(X)\) maps an
  object \(D\) of \(\DivQnn(X,L)\) to the support of the
  multiplier ideal sheaf for \(D\).
  \begin{equation}
    \label{eqn:ample_homotopy:3}
    \sMI(D)=\supp\MI(D)
  \end{equation}
  I.e., for each \(x\in X\), \(x\in \sMI(D)\) iff
  \(\MI(D)_{x}= \O_{X,x}\).  Then, the following are
  equivalent.
  \begin{enumerate}
  \item \(L\) is ample.
  \item\label{item:6} \(L\) satisfies (SoO).
  \item\label{item:7} \(\sMI^{int}\) is surjective on objects.
  \item\label{item:1} \(\sMI^{int}\) is essentially surjective on
    objects.
  \end{enumerate}
  (1) and (2) are equivalent by
  Theorem~\ref{thm:ample_homotopy:1}. (2) and (3) are
  equivalent by the definition
  \eqref{eqn:presheaves_and_Q-divisors:82}.  (3) and (4) are
  equivalent because every isomorphism in \(\Opz(X)\) is an
  identity.  So Theorem~\ref{thm:ample_homotopically:1}
  follows from Lemma~\ref{lem:presheaves_and_Q-divisors:15},
  Lemma~\ref{lem:presheaves_and_Q-divisors:16} and
  Theorem~\ref{thm:presheaves_and_Q-divisors:6}.
\end{proof}

\begin{remark}
  \label{rem:ample_homotopy:5}
  In the proof of Theorem~\ref{thm:ample_homotopically:1},
  we have made a choice for the functor
  \(\sMI:\DivQnn(X,L)\rightarrow \Opz(X)\). But it is not
  unique or canonical. We can also use the non-lc ideal
  sheaves instead (\cite{fujino-10},
  \cite{fujino-schwede-takagi-11}). It might be possible to
  use a totally different invariants of \(\bbQ\)-divisor to
  define the functor \(\sMI:\DivQnn(X,L)\rightarrow \Opz(X)\)
  in such a way that the smoothness can be removed from the
  assumption.
\end{remark}


\section{Examples of admissible subcategories}
\label{sec:example-admiss-subc}

Let \(X\) be a complex smooth projective variety. Let \(L\)
be an ample divisor. Then the functor \(\sMI\) in
Theorem~\ref{thm:ample_homotopically:1} generates a
homotopical presentation of \(\Opz(X)\). Here the relation
is the set of all morphisms in \(\MM(X,L)\) mapped to an
isomorphism in \(\Opz(X)\). We can make it smaller using
Theorem~\ref{thm:presheaves_and_Q-divisors:1}. Instead of
repeating with \(H\), we demonstrate it using the
restriction \(\MM(X,L,x)\rightarrow \Opz(X,x)\) of \(\sMI\)
in \eqref{eqn:ample_homotopy:25} below.

Let \(x\) be a point of \(X\). We denote by
\begin{equation}
  \label{eqn:ample_homotopy:26}
  \DivQnn(X,L,x)
\end{equation}
the full subcategory of \(\DivQnn(X,L)\) such that for an
object \(D\) of \(\DivQnn(X,L)\), \(D\) is an object of
\(\DivQnn(X,L,x)\) iff \(x\not\in \sMI(D)\). We denote by
\begin{equation}
  \label{eqn:ample_homotopy:27}
  \MM(X,L,x)
\end{equation}
the full subcategory of \(\MM(X,L)\) such that for an object
\(\uD\) of \(\MM(X,L)\), \(\uD\) is an object of
\(\MM(X,L,x)\) iff every element of \(\uD\) is an object of
\(\DivQnn(X,L,x)\).  We denote by
\begin{equation}
  \label{eqn:ample_homotopy:28}
  \Opz(X,x)
\end{equation}
the full subcategory of \(\Opz(X)\) on open sets \(U\) with
\(x\not\in U\).  Both of \(\MM(X,L,x)\) and \(\Opz(X,x)\)
are closed under finite limits.

We use the same notation \(\sMI\) for the restrictions
\begin{equation}
  \label{eqn:ample_homotopy:25}
  \sMI:\MM(X,L,x)\rightarrow \Opz(X,x)
\end{equation}
and 
\begin{equation}
  \label{eqn:ample_homotopy:29}
  \sMI:\DivQnn(X,L,x)\rightarrow \Opz(X,x).
\end{equation}
For \(D\in\DivQnn(X,L,x)\),
denote by
\begin{equation}
  \label{eqn:ample_homotopy:30}
  c_{x}(D)
\end{equation}
the log-canonical threshold for \(D\) at \(x\)
(Definition~9.3.12 in~\cite{lazarsfeld-04}).  If \((X,D)\)
is log canonical, we denote by
\begin{equation}
  \label{eqn:ample_homotopy:31}
  \mc_{x}(D)
\end{equation}
the minimal center of \(D\) at \(x\) (\cite{kawamata-97}).

\begin{definition}
  \label{def:ample_homotopy:8}
  Let \(X\) be a complex smooth projective variety. Let
  \(L\) be an ample divisor on \(X\).  Let \(x\) be a point
  of \(X\).
  \begin{enumerate}
  \item We call a morphism \(\phi:E\rightarrow D\) in
    \(\Div_{\bbQ}(X,L,x)\) \textbf{reducing} if the
    following hold.
    \begin{enumerate}
    \item\label{item:12} \(\sMI(\phi)\) is an identity
      morphism in \(\Opz(X,x)\).
      \begin{equation}
        \sMI(D)=\sMI(E)
      \end{equation}        
    \item\label{item:13} If the dimension of the minimal
      center of \(c_{x}(D)\cdot D\) at \(x\) is positive,
      \begin{equation}
        \mc_{x}(c_{x}(E)\cdot E)\subsetneqq \mc_{x}(c_{x}(D)\cdot D).
      \end{equation}
      Otherwise,
      \begin{equation}
        \label{eqn:ample_homotopy:32}
        \mc_{x}(c_{x}(E)\cdot E)=\{x\}= \mc_{x}(c_{x}(D)\cdot D).
      \end{equation}
     \end{enumerate}
   \item We call a morphism \(\phi:\uE\rightarrow \uD\) in
     \(\MM(X,L,x)\) \textbf{reducing} if the following hold.
    \begin{enumerate}
    \item\label{item:2} \(\sMI(\phi)\) is the identity
      morphism in \(\Opz(X,x)\), i.e.,
      \(\phi\in\MM(X,L,x)_{\uD}\)
      (Definition~\ref{def:homotopical_presentation:2}).
    \item\label{item:3} For each \(D_{i}\in \uD\) there is
      \(E_{j}\in \uE\) and a reducing morphism
      \(E_{j}\rightarrow D_{i}\).
    \end{enumerate}
  \end{enumerate}
\end{definition}

\begin{definition}
  Let \(X\) be a complex smooth projective variety. Let
  \(L\) be an ample divisor on \(X\).  Let \(x\) be a point
  of \(X\).  For each \(\uD\in\MM(X,L,x)\), we denote by
  \begin{equation}
    \label{eqn:ample_homotopy:33}
    A(X,L,x)_{\uD}
  \end{equation}
  the full subcategory of \(\MM(X,L,x)_{\uD}\) on the pairs
  \(\gangle{\uE,\phi}\) of \(\uE\in\MM(X,L,x)\) and a
  reducing morphism \(\phi:\uE\rightarrow \uD\).
\end{definition}

\begin{definition}
  Let \(X\) be a complex smooth projective variety. Let
  \(L\) be an ample divisor on \(X\). Let \(x\) be a point
  of \(X\). Let \(B\in\DivQnn(X,L)\) and
  \(D\in\DivQnn(X,L,x)\).  We say that \(B\) is a
  \textbf{reducing with respect to \(D\) at \(x\)} if the
  morphism
  \begin{equation}
    \label{eqn:ample_homotopy:18}
    (B+D)\rightarrow D
  \end{equation}
  in \(\DivQnn(X,L,x)\) is reducing.
\end{definition}

The following lemma is well-known (cf.~Proposition~10.4.10
in~\cite{lazarsfeld-04}).

\begin{lemma}
  \label{lem:ample_homotopy:4}
  Let \(X\) be a complex smooth projective variety. Let
  \(L\) be an ample divisor on \(X\). Let \(x\) be a point
  of \(X\). Let \(D\in\DivQnn(X,L)\). Assume that
  \begin{equation}
    \dim\mc_{x}(c_{x}(D)\cdot D)>0.
  \end{equation}  
  Then there exists \(B\in \DivQnn(X,L)\) reducing with
  respect to \(D\) at \(x\).
\end{lemma}

\begin{lemma}
  \label{lem:ample_homotopy:5}
  Let \(X\) be a complex smooth projective variety. Let
  \(L\) be an ample divisor on \(X\).  Let \(x\) be a point
  of \(X\).  For each \(\uD\in\MM(X,L,x)\),
  \(A(X,L,x)_{\uD}\) is admissible
  (Definition~\ref{def:ample_homotopy:3}). I.e., the
  following hold.
  \begin{enumerate}
  \item \(A(X,L,x)_{\uD}\) is cofiltered.
  \item \(\sMI\) in \eqref{eqn:ample_homotopy:25}
    has \(A(X,L,x)_{\uD}\)-lifting property
    (Definition~\ref{def:homotopical_presentation:2}).
  \end{enumerate}
\end{lemma}
\begin{proof}
  (1) Since \(\MM(X,L,x)\) is a preorder, it is enough to
  show that it is closed under products. Let
  \((\uE_{1},\phi_{1})\) and \((\uE_{2},\phi_{2})\) be two
  objects of \(A(X,L,x)_{\uD}\). Write
  \(\uE=\uE_{1}\cup\uE_{2}\). Let
  \(\phi:\uE\rightarrow \uD\) be the canonical morphism so
  that \((\uE,\phi)\) is the product of
  \((\uE_{1},\phi_{1})\) and \((\uE_{2},\phi_{2})\).  Then
  \(\sMI(\uE_{1})=\sMI(\uD)\) and
  \(\sMI(\uE_{2})=\sMI(\uD)\) imply
  \(\sMI(\uE)=\sMI(\uE_{1}\cup \uE_{2})=\sMI(\uE_{1})\cap
  \sMI(\uE_{2})= \sMI(\uD)\). Thus \eqref{item:2} in
  Definition~\ref{def:ample_homotopy:8} holds for \(\phi\).
  \eqref{item:3} holds by definition.

  (2) Let \(\uE\) be an object of \(\MM(X,L,x)\). Let
  \(\delta:\sMI(\uD)\rightarrow \sMI(\uE)\) be a
  morphism in \(\Opz(X,x)\), i.e., we have an inclusion
  \begin{equation}
    \sMI(\uD)\subseteq \sMI(\uE).
  \end{equation}
  Let \(\uD=\{D_{1},\dots,D_{n}\}\). We define an object
  \(\uE'=\{E_{1},\dots,E_{n}\}\) of \(\MM(X,L,x)\) by
  \begin{equation}
    E_{i}=D_{i}+B_{i}
  \end{equation}
  where \(B_{i}\) is reducing with respect to \(D_{i}\) at
  \(x\) (Lemma~\ref{lem:ample_homotopy:4}) so that the
  canonical morphisms \(E_{i}\rightarrow D_{i}\) is reducing
  for all \(i\).

  Let \(\uE''=\uE'\cup \uE\).  The canonical morphisms
  \(\phi:\uE''\rightarrow \uD\) and
  \(\psi:\uE''\rightarrow \uE\) satisfy
  \(\sMI(\psi)=\delta\cdot \sMI(\phi)\). The morphism
  \(\phi\) satisfies \eqref{item:2} and \eqref{item:3} by
  construction. Thus
  \(\gangle{\uE,\phi}\in A(X,L,x)_{\uD}\). Therefore
  \(\sMI\) has \(A(X,L,x)_{\uD}\)-lifting property.
\end{proof}

\begin{proposition}
  \label{pro:ample_homotopy:1}
  Let \(X\) be a complex smooth projective variety. Let
  \(L\) be an ample divisor on \(X\). Let \(x\) be a point
  of \(X\).  Then the Quillen adjunction
  \begin{equation}
    \label{eqn:ample_homotopy:34}
    \sMI_{*}:\BL_{S_{A(X,L,x)}}(\Un{\MM(X,L,x)})\rightleftarrows
    \Un{\Opz(X,x)}:\sMI^{*}
  \end{equation}
  associated with \(\sMI\) in \eqref{eqn:ample_homotopy:25}
  is a Quillen equivalence.
\end{proposition}

\providecommand{\bysame}{\leavevmode\hbox to3em{\hrulefill}\thinspace}
\providecommand{\MR}{\relax\ifhmode\unskip\space\fi MR }
\providecommand{\MRhref}[2]{%
  \href{http://www.ams.org/mathscinet-getitem?mr=#1}{#2}
}
\providecommand{\href}[2]{#2}

\end{document}